%Paper: funct-an/9211003
%From: arveson@math.berkeley.edu (Bill Arveson)
%Date: Sun, 22 Nov 92 16:00:36 PST
%Date (revised): Tue, 24 Nov 92 07:19:04 PST

%
%NOTE: This document requires AMSTeX 2.0 or 2.1, the
%accompanying version of the amsppt style macros,
%and the font package amsfonts 2.0 or 2.1. If you try
%to typeset it using an older version of AMSTeX, it
%probably won't work.  You can probably get by
%without the fonts package, but you will surely need the
%font metrics for amsfonts 2.0 or 2.1, preferably 2.1.
%
%NOTE: The entire AMSTeX 2.1 package can be obtained
%(free) from the AMS's ftp site in Providence. If you
%are running AMSTeX 2.0, the AMS strongly recommends that
%you upgrade to 2.1.  In particular, this will resolve some
%difficulties that 2.0 has with some fonts, and will
%improve backward compatibility with older versions of TeX.
%
%
%This paper is a sequel to "C*-algebras and numerical
%linear algebra".  It has been submitted to the Acta
%Sci Math (Szeged).
%
%%%%%%%%%%%%%%%%%%%  BEGINNING OF TEXFILE  %%%%%%%%%%%
%
\input amstex
%Comment the preceding line out if your system loads the
%AMSTeX format automatically.
%
%
%MACROS
\def\cstar{$C^*$-algebra}

\def\I{\bold1}

\def\deg{\roman{deg }}
\def\rank{\roman{rank}}
%
%AMSTeX HEADERS
\documentstyle{amsppt}
%Comment the preceding line out if your system loads
%the amsppt style macros automatically with the amstex format.
\loadbold
\magnification=\magstep 1
\topmatter
\title
Improper filtrations for $C^*$-algebras:\\
spectra of unilateral tridiagonal operators
\endtitle
\rightheadtext{Improper filtrations for $C^*$-algebras}
\author William Arveson
\endauthor
\email
arveson\@math.berkeley.edu
\endemail
%
%\date
%11 February 1992
%\enddate
%
\affil
Department of Mathematics\\
University of California\\
Berkeley, CA 94720 USA
\endaffil
\dedicatory
To B\'ela Sz\"okefalvi-Nagy on the occasion of his eightieth birthday
\enddedicatory
\thanks
This research was supported in part by
NSF grants DMS89-12362 and DMS-9212893
\endthanks
\abstract Let $\Cal A\subseteq\Cal B(H)$ be a \cstar\ of operators
and let $P_1\leq P_2\leq\dots$ be an increasing sequence of
finite dimensional projections in $\Cal B(H)$.
In a previous paper \cite{3} we developed methods for computing
the spectrum of self adjoint operators $T\in\Cal A$ in terms of
the spectra of the associated sequence of finite dimensional
compressions $P_nTP_n$.  In a suitable context, we
showed that this is possible when $P_n$ increases to $\I$.  In this
paper we drop that hypothesis and obtain an appropriate
generalization of the main results of \cite{3}.

Let $P_+=\lim_nP_n$, $H_+=P_+H$.  The set $\Cal A_+\subseteq\Cal B(H_+)$
of all compact perturbations of operators $P_+T\!\restriction_{H_+}$,
$T\in\Cal A$, is a
\cstar\ which is somewhat analogous to the Toeplitz \cstar\ acting
on $H^2$.
Indeed, in the most important examples $\Cal A$ is a simple unital
\cstar\ having a unique tracial state, the operators in $\Cal A$ are
``bilateral", those in $\Cal A_+$ are ``unilateral", and there is a
short exact sequence of \cstar s
$$
0\to\Cal K\to\Cal A_+\to\Cal A\to 0
$$
whose features are central to this problem of approximating spectra
of operators in $\Cal A$ in terms of the eigenvalues of their finite
dimensional compressions along the given filtration.

This work was undertaken in order to develop an efficient method
for computing the spectra of discretized Hamiltonians of one
dimensional quantum
systems in terms of ``unilateral" tridiagonal $n\times n$ matrices.
The solution of that problem is presented in Theorem 3.4.
\endabstract
\endtopmatter
\vfill\eject
%Replace \pagebreak below with the line above
%to fill the lower part of the title page with
%space, rather than stretching it.
\pagebreak
\document
\subheading{1. Introduction}

There is a natural way to discretize the Hamiltonian of a
one dimensional quantum system in such a way that
a) the uncertainty principle is preserved, and b) the
discretized operator is in a form appropriate for carrying out
numerical studies \cite{1,2}.  In a suitable representation,
this discretized Hamiltonian is a bilateral tridiagonal
operator $T$ of the form
$$
Te_n = e_{n-1} + v(cos(n\theta))e_n + e_{n+1} \leqno{1.1}
$$
where $\{e_n: n\in \Bbb Z\}$ is an orthonormal basis
for a Hilbert space $H$, and $v$ is a real function
in $C[-1,+1]$.  $v$ is a rescaled version of
the potential of the system.  The number $\theta$ is
related to the numerical step size, and may properly be considered
an irrational multiple of $\pi$ \cite{1}.

The spectra of these operators cannot be computed explicitly.
For example, even in the simplest case where $v(x) = 2x$ (in physical
terms, the case
of a one dimensional harmonic oscillator) it is not yet
known if the spectrum of $T$ is totally disconnected.
Choi, Elliott and Lui have shown that this is the case
when $\theta/\pi$ is a Liouville number \cite{9}, but despite
significant progress since the early eighties this
`Ten Martini' problem of Mark Kac remains open
in general \cite{5,6,7,8,9,14}.  To our knowledge,
essentially nothing is known in the case of an arbitrary
continuous potential $v$.  Thus it is of
interest to understand how one might carry
out numerical calculations of the spectrum using
established techniques for computing eigenvalues of self-adjoint
$n\times n$ matrices, and to describe the precise sense in
which these eigenvalues converge to the spectrum of $T$ as $n\to\infty$.
A general approach to this problem was worked out in \cite{3}.

Here we will consider operators $T$ of the following
somewhat more general tridiagonal form
$$
Te_n = e_{n-1} + d_ne_n + e_{n+1} \leqno{1.2}
$$
where $\{d_n: n\in \Bbb Z\}$ is a bounded bilateral sequence
of real numbers which is {\sl almost periodic but not periodic}.
A sequence $d_n = v(cos(n\theta))$ as in formula 1.1 will
satisfy this hypothesis whenever $\theta/\pi$ is irrational
and $v\in C[-1,+1]$ is not a constant.
In \cite{3}, we showed that
such a program will succeed for operators of the form 1.1
by establishing a formula which
shows how the spectrum of $T$ is determined by the
eigenvalues of the sequence of $2n+1\times 2n+1$ matrices
$$
\pmatrix
d_{-n}&1&0&\hdots&0&0\\
1&d_{-n+1}&1&\hdots&0&0\\
0&1&d_{-n+2}&\hdots&0&0\\
\vdots&\vdots&\vdots&\ddots&\vdots&\vdots\\
0&0&0&\hdots&d_{n-1}&1\\
0&0&0&\hdots&1&d_n
\endpmatrix.
$$
However, it is more convenient to compute with the sequence
of ``unilateral" $n\times n$ compressions of $T$
$$
T_n =
\pmatrix
d_{1}&1&0&\hdots&0&0\\
1&d_{2}&1&\hdots&0&0\\
0&1&d_{3}&\hdots&0&0\\
\vdots&\vdots&\vdots&\ddots&\vdots&\vdots\\
0&0&0&\hdots&d_{n-1}&1\\
0&0&0&\hdots&1&d_n
\endpmatrix
{}.
\leqno{1.3}
$$
In fact, the latter has been implemented for the Ten Martini
operators of Mark Kac in a Macintosh program \cite{4}, and it
is effective.  In this paper we will prove that this
procedure will {\sl always} be successful by establishing
a variation of the formula \cite{3, formula 5.1}
which is applicable to the sequence
of unilateral compressions $1.3$.

In this setting the formula takes the following form.  Let
$T_n$ be as in $1.3$, $n = 1, 2, \dots$ and let
$\lambda_1^n < \lambda_2^n <\dots < \lambda_n^n$ be the eigenvalue
list of $T_n$ (recall that the eigenvalues are distinct because
$T_n$ is tridiagonal \cite{12, p. 124}).  Then there is a probability
measure $\mu_T$, whose closed support is precisely the spectrum
of $T$, such that for every continuous function $f \in C_0(\Bbb R)$
we have
$$
\lim_{n\to\infty}{1\over n}(f(\lambda_1^n) +
f(\lambda_2^n) +\cdots +f(\lambda_n^n)) =
\int_{-\infty}^{+\infty} f(x)\,d\mu_T(x).  \leqno{1.4}
$$
It follows that a real number $\lambda$ is in the spectrum of $T$
iff the eigenvalues of the sequence of matrices $T_1, T_2, \dots$
accumulate rapidly enough in any fixed neighborhood $U$ of $\lambda$
so that their density in $U$ is positive.  We also describe the measure
$\mu_T$ in terms of of the unique tracial state on a certain
\cstar\ associated with the operator $T$.  See Theorem 3.4 and
the remarks following it.

Formula $1.4$ shows quite explicitly how one should go about
calculating approximations to the measure $\mu_T$, and it is
these calculations that have been implemented in the
program \cite{4} cited above.
\subheading{2. Improper filtrations}

In this section we analyze the situation described above in a
rather abstract context, and establish an appropriate generalization
of the main result of \cite{3}.
Let us first recall some terminology from \cite{3}.  Let $H$
be a separable infinite dimensional Hilbert space.
A {\sl filtration} of $H$ is an increasing sequence
$H_1\subseteq H_2\subseteq\dots$ of finite dimensional subspaces of
$H$ with the property that $\cup_n H_n$ is dense in $H$.

Here, we will simply drop the latter condition on the sequence
$H_n$.  In order to avoid annoying trivialities, we do require
that the closed subspace
$$
H_+ = \overline{\cup_n H_n}
$$
generated by the filtration should be infinite dimensional.
The filtration is called {\sl proper} or {\sl improper}
according as $H_+ = H$ or $H_+ \neq H$.  Let $P_n$ denote the
projection onto $H_n$.  As in \cite{3} we may introduce the
notion of the {\sl degree} of an operator $T\in\Cal B(H)$ relative
to a filtration,
$$
\deg(T) = \sup_n \rank(P_nT - TP_n).
$$
The degree of an operator is a nonnegative integer
or $+\infty$, and properties
established in \cite{3} such as $\deg(ST)\leq \deg(S) + \deg(T)$
persist in this more general context.

For example, if $\{e_1, e_2, \dots\}$ is an infinite orthonormal
sequence in $H$ and $H_n$ is the linear span of $\{e_1,\dots, e_n\}$,
then $\{H_1\subseteq H_2\subseteq\dots\}$ is a filtration which is typically
improper.  As we pointed out in \cite{3} in the case of proper
filtrations, any operator whose matrix
relative to this orthonormal set is band-limited must have finite degree.
Thus, {\sl it is as appropriate here as it was in \cite{3}
to think of finite degree operators as abstractions
of band-limited matrices}.  Any self-adjoint operator
$T$ determines a multitude of filtrations relative to which it
has degree 1: choose any vector $e$
in $H$ such that $\{e, Te, T^2e, \dots\}$ is linearly independent
and put $H_n = [e, Te, T^2e,\dots, T^{n-1}e]$.  This filtration
is proper iff $e$ is a cyclic vector for $T$.  Of course, in general
there can be very irregular jumps in the dimensions of the subspaces of
a filtration.

Finally, if $\Cal A$ is a \cstar\ of operators acting on $H$
and $\Cal F = \{H_1 \subseteq H_2 \subseteq \dots\}$ is a filtration
of $H$, then $\Cal F$ is called an $\Cal A${\sl-filtration} if the
set of all finite degree operators in $\Cal A$ is norm-dense in
$\Cal A$.  If $\Cal A$ is generated by a set of finite degree
operators, then $\Cal F$ is called an $\Cal A$-filtration \cite{3}.
We begin with two propositions which concern compact perturbations.

\proclaim{Proposition 2.1}Let $\Cal A\subseteq\Cal B(H)$ be
a \cstar, let $\{H_1\subseteq H_2\subseteq\dots\}$ be an
$\Cal A$-filtration and let $P_+$ be the projection onto
$H_+$.  Then $P_+A - AP_+$ is compact for every $A\in\Cal A$.
\endproclaim
\demo{proof}  Since the finite degree operators in $\Cal A$ are
norm-dense and $P_+A-AP_+$ is norm-continuous in $A$,
it suffices to show that
$$
\rank(P_+A - AP_+) \leq \deg(A) <\infty
$$
for every $A\in\Cal A$ having finite degree.  Since
$$
\rank(P_nA - AP_n) \leq \deg(A)
$$
for every $n = 1, 2,\dots$ and since $P_nA - AP_n$ converges
to $P_+A - AP_+$ in the strong operator topology, the propostion
follows from the fact that the rank function is weakly lower
semicontinuous on $\Cal B(H)$.

Since we lack an appropriate reference for
the latter, we sketch a short proof for completeness.
The assertion is that
for every positive integer $k$, the set of operators
$$
\Cal S_k = \{T\in\Cal B(H): \rank(T)\leq k\}
$$
is weakly closed.  Let $K$ be the Hilbert space
antisymmetric tensor product
of $k+1$ copies of $H$.  A subspace $M$ of $H$ has dimension
at most $k$ iff
$$
M^{\wedge k+1} = M\wedge M\wedge\dots\wedge M = 0.
$$
Thus, an operator $T$ belongs to $\Cal S_k$ iff
$$
T\xi_1\wedge T\xi_2\wedge\dots\wedge T\xi_{k+1} = 0
$$
for every choice of vectors $\xi_1,\xi_1,\dots\xi_{k+1}$ in $H$.
Using the fact that $K$ is spanned by vectors of the form
$\eta_1\wedge\eta_2\wedge\dots\wedge\eta_{k+1}$
we see that $T$ belongs to $\Cal S_k$
iff for every set of $2k+2$ vectors
$\xi_1,\dots,\xi_{k+1},\eta_1,\dots\eta_{k+1}$ in
$H$, the polynomial $p(T)$ defined by
$$
p(T) = <T\xi_1\wedge T\xi_2\wedge\dots\wedge T\xi_{k+1},
\eta_1\wedge\eta_2\wedge\dots\wedge\eta_{k+1}>
$$
is zero.  Taking note of the inner product in $K$,
we see that $p(T)$ is the determinant
of the $k+1\times k+1$
matrix whose $ij$th term is $<T\xi_i,\eta_j>$, and is
therefore weakly continuous.  Hence $\Cal S_k$ is weakly closed \qed
\enddemo

If one is given a positive linear functional $\rho$ on a
\cstar\ $\Cal C$, then the GNS construction gives rise to
a representation $\pi_\rho$ of $\Cal C$ on a Hilbert space
$H_\rho$.  We will say that $\rho$ is {\sl infinite
dimensional} if $H_\rho$ is an infinite dimensional Hilbert
space.

\proclaim{Proposition 2.2}Let $\Cal A\subseteq\Cal B(H)$ be
a unital \cstar\ having a unique tracial state $\tau$, and
assume $\tau$ is infinite dimensional.  Then $\tau(K) = 0$
for every compact operator  $K\in\Cal A$.
\endproclaim
\demo{proof}  Let $\pi$ be the representation of $\Cal A$
obtained via the GNS construction on $\tau$.
Then $M = \pi(\Cal A)^{\prime\prime}$
is a von Neumann algebra with the property that the natural extension
$\tilde\tau$ of $\tau$ to $M$ is a normal tracial state of $M$.

Notice that $M$ is a factor.  Indeed, if $e$ is any nonzero
central projection in $M$ then note first that $t=\tilde\tau(e)>0$.
Indeed, by the Schwarz inequality, $\tilde\tau(e)=0$ implies
$\tilde\tau(aeb)=\tau(eba)=0$ for every $a,b\in M$,
and hence $e=0$ because
$\tilde\tau$ is a vector state of $M$ which is associated
with a cyclic vector.  In order to
show that $e=\I$, we define a tracial state $\rho$ of $\Cal A$ by
$$
\rho(A) = \frac{1}{t}\tilde\tau(\pi(A)e).
$$
Because of uniqueness, we have $\rho = \tau$.
{}From this, together with the fact that $\pi(\Cal A)$ is strongly
dense in $M$, it follows that
$\tilde\tau(ae) = t\tilde\tau(a)$ for every
$a\in M$.  By an argument similar to the one given above to show
that $e\neq 0$, we find that $e-t\I=0$, hence $e=t\I$.
Thus $t=1$ and $e=\I$.

We may conclude that $M$ is a finite factor.
Now let $\Cal J$ be the ideal of all compact operators in $\Cal A$.
$\Cal J$ is a \cstar\ of compact operators, and hence it must
decompose uniquely into a \cstar ic direct sum of elementary \cstar s
$$
\Cal J = \sum_k\Cal J_k.
$$
It suffices to show that $\tau(\Cal J_k)=0$ for every $k$.
Contrapositively,
assume $\tau(\Cal J_k)\neq 0$ for some $k$.  Then $\pi(\Cal J_k)$
must be nonzero too, and its weak closure in $M$ must be all of
$M$ since it is a nonzero
weakly closed ideal in a factor.  Because $\tau$ is an infinite
dimensional state $M$ is infinite dimensional and therefore so is
$\Cal J_k$.  But an infinite-dimensional elementary \cstar\ has no
nonzero bounded traces.  Hence the restriction of $\tau$ to $\Cal J_k$
must be zero \qed
\enddemo

\remark{Remark} Perhaps it is worth pointing out that Proposition 2.2
becomes false if one deletes the hypothesis that $\tau$ is an infinite
dimensional state.  For example, choose any \cstar\ of the form
$A = B\oplus C$ where $B$ is a unital \cstar\ having no nonzero
finite traces (such as a Cuntz algebra) and $C$ is an $n\times n$
matrix algebra.  Certainly $A$ has a unique tracial state $\tau$,
and there is an obvious faithful representation of $A$ which
contains finite rank operators not annihilated by $\tau$.
\endremark

Throughout the remainder of this section, we assume that we are
given a unital \cstar\ $\Cal A$ which has a unique tracial state
$\tau$.  Every self-adjoint operator $A\in\Cal A$ gives rise to
a probability measure $\mu_A$ on the real line by way of
$$
\int_{-\infty}^{+\infty}f(x)\,d\mu_A(x) = \tau(f(A)),\qquad f\in C_0(\Bbb R).
$$
$\mu_A$ is called the {\sl spectral distribution} of $A$
\cite{3, section 4}.  If $\tau$ is faithful then the closed support of
$\mu_A$ is precisely the spectrum of $A$.
Let $\Cal F = \{H_1\subseteq H_2\dots\}$ be a perhaps improper
$\Cal A$-filtration and let $H_+$ and $P_+$ have the
meaning assigned above.  The following result is a
generalization of theorem 4.5 of \cite{3} to the case of improper
filtrations.

\proclaim{Theorem 2.3} Let $\Cal A\subseteq\Cal B(H)$ be a
\cstar\ having a unique tracial state, let
$\{H_1\subseteq H_2\subseteq\dots\}$ be an $\Cal A$-filtration,
and assume that the
subspace $H_+$ has the following property
$$
A\!\restriction_{H_+} = {\roman compact}\implies A=
{\roman compact},\leqno{2.4}
$$
for every $A\in\Cal A$.  Let $\mu_A$ be the
spectral distribution of a self-adjoint operator $A\in\Cal A$, and
let $[a,b]$ be the smallest closed interval containing $\sigma(A)$.
For each $n$, let $d_n=\dim H_n$ and
let $\lambda_1^n,\lambda_2^n,\dots,\lambda_{d_n}^n$ be the eigenvalue list of
$A_n=P_nA\!\restriction_{H_n}$.  Then for every $f\in C[a,b]$,
$$
\lim_{n\to\infty}{1\over
{d_n}}(f(\lambda_1^n)+f(\lambda_2^n)+\dots+f(\lambda_{d_n}^n))
= \int_a^b f(x)\, d\mu_A(x).
$$
\endproclaim
\remark{Remarks}
We will write $\Cal K$ and
$\Cal K_+$ for the respective algebras of compact operators on
$H$ and $H_+$.  Notice that the hypothesis 2.4 is satisfied for
trivial reasons if the filtration is proper, and in that case
Theorem 2.3 reduces to \cite{3, Theorem 4.5}.  2.4 is also automatic
when $\Cal A$ is simple. Indeed, if $T$ is an operator in $\Cal A$
whose restriction to $H_+$ is compact, then $T$ belongs to the kernel
of the $*$-homomorphism $\omega :\Cal A\to\Cal B(H_+)/\Cal K_+$ given
by $\omega(X) = P_+X\!\restriction_{H_+} + \ \Cal K_+$.  Since $\omega$
is not the zero homomorphism ($H_+$ is infinite dimensional),
it must be injective by simplicity.  Hence $T=0$.
\endremark

\demo{proof of Theorem 2.3}
Let $\Cal A_+$ be the set of operators in $\Cal B(H_+)$ of the
form
$$
\Cal A_+=P_+A\!\restriction_{H_+}+\ K,
$$
where $A\in\Cal A$ and $K$ is
a compact operator on $H_+$.  We will show first that

\roster
\item"{2.5}" $\Cal A_+$ is a \cstar\ algebra having a unique tracial
state $\tau_+$.
\item"{2.6}" $\tau_+$ is related to the tracial state $\tau$ of
$\Cal A$ by way of
$$
\tau_+(P_+A\!\restriction_{H_+}+\ K) = \tau(A),\qquad A\in\Cal A, K\in\Cal K_+.
$$
\endroster

If $A$ and $B$ are two
operators in $\Cal A$ with $A_+$ and $B_+$ denoting their respective
compressions to $H_+$, then from Proposition 2.1 we see that
$$
A_+B_+ = P_{H_+}AB\!\restriction_{H_+} + \ \roman{compact}.
$$
It follows that $\Cal A_+$ is a unital $*$-subalgebra of $\Cal B(H_+)$
which contains all compact operators.  Thus, to prove that $\Cal A_+$
is a \cstar, it suffices to show that its image in the Calkin algebra
$\Cal B(H_+)/\Cal K_+$ of $H_+$ is norm-closed.  But the map
$\pi:\Cal A\to \Cal B(H_+)/\Cal K_+$ defined by
$$
\pi(A) = P_{H_+}A\!\restriction_{H_+} +\ \Cal K_+
$$
is clearly a unital $*$-homomorphism
whose range is the image of
$\Cal A_+$ in $\Cal B(H_+)/\Cal K_+$.  Since the range of $\pi$
must be norm-closed, we conclude that $\Cal A_+$ is norm-closed.

In order to prove 2.5, we first exhibit a tracial
state of $\Cal A_+$.  Let $\sigma:\Cal A\to\Cal B(H_+)/\Cal K_+$
be the map defined by
$$
\sigma(A) = P_+A\!\restriction_{H_+} +\ \Cal K_+.
$$
By Proposition 2.1, $\sigma$ is a $*$-homomorphism of unital
\cstar s, and its range is precisely the image of $\Cal A_+$ in
the Calkin algebra.  By 2.4, the kernel of $\sigma$ is
$\Cal A\cap\Cal K$; indeed $\sigma(A)=0$ iff $\sigma(A^*A)=0$
iff the restriction of $A$ to $H_+$ is compact iff $A$ is
compact.  Hence $\sigma$ determines an isomorphism
$\dot\sigma$ of \cstar s
$$
\dot\sigma: \Cal A/\Cal A\cap\Cal K \to \Cal A_+/\Cal K_+.
$$
On the other hand, because of Proposition 2.2, $\tau$ vanishes
on $\Cal A\cap\Cal K$ and thus can be promoted to a tracial state
$\dot\tau$ on the quotient $\Cal A/\Cal A\cap\Cal K$ by way of
$\dot\tau(A + \Cal A\cap\Cal K)=\tau(A)$.  Thus the composition
$\dot\tau\circ\dot\sigma^{-1}$ defines a tracial state of
$\Cal A_+/\Cal K_+$.  If we now define $\tau_+$ on $\Cal A_+$ by
$\tau_+(A) = \dot\tau\circ\dot\sigma^{-1}(A+\Cal K_+)$, then $\tau_+$
is certainly a tracial state on $\Cal A_+$.
After unravelling its definition we find that
$\tau$ is related to $\tau_+$ by formula 2.6
$$
\tau_+(P_+A\!\restriction +\ \Cal K_+) = \tau(A),\qquad A\in\Cal A.
\leqno{2.7}
$$

It remains to show that $\tau_+$ is the only tracial state on $\Cal A_+$.
To see that, let $\rho$ be a tracial state on $\Cal A_+$.
It is clear that the restriction of $\rho$ to
$\Cal K_+$ must be zero, because there are no nonzero bounded traces
on the algebra of compact operators on an infinite dimensional
Hilbert space.  Hence we may promote $\rho$ to a tracial state
$\dot\rho$ on $\Cal A_+/\Cal K_+$ by
$$
\dot\rho(T+\Cal K_+) = \rho(T),\qquad T\in\Cal A_+.
$$
Making use of the isomorphism
$\sigma:\Cal A/{\Cal A\cap\Cal K}\to\Cal A_+/\Cal K_+$ exhibited above,
we define a state $\rho_0$ on $A$ by
$$
\rho_0(A) = \dot\rho\circ\sigma(A+\Cal A\cap\Cal K)=
\rho(P_+A\!\restriction_{H_+} +\ \Cal K_+).
\leqno{2.8}
$$
$\rho_0$ is clearly a trace on $A$, hence $\rho_0=\tau$.  2.7 and 2.8
together show that $\rho$ is the trace $\tau_+$ exhibited above.  Hence
both assertions 2.5 and 2.6 are established.

It remains to deduce Theorem 2.3 from these two assertions and the
results of \cite{3}.  For that, we view
$\Cal F_+=\{H_1\subseteq H_2\subseteq\dots\}$ as a
{\sl proper} filtration of $H_+$.  Notice that if $A$ is any finite
degree operator in $\Cal A$ and $B$ is an operator in $\Cal A_+$
of the form $B=P_+A\!\restriction_{H_+} + F$ where $F$ is a finite
rank operator on $H_+$, then for every $n = 1,2,\dots$ we have
$P_nB-BP_n=P_+(P_nA-AP_n)\!\restriction_{H_+} + P_nF - FP_n$.  Hence
$$
\rank(P_nB-BP_n)\leq \sup_n rank(P_nA-AP_n) + 2\rank(F)
\leq\deg(A) + 2\rank(F)<\infty.
$$
It follows that the degree of $B$ is finite.  Such operators
$B$ are norm-dense in $\Cal A_+$ since the finite degree operators
of $\Cal A$ are norm dense in $\Cal A$.  We conclude that $\Cal F$ {\sl is
a proper} $\Cal A$--{\sl filtration of} $H_+$.

Thus we are in position to apply \cite{3, Theorem 4.5} to the operator
$A_+=P_+A\!\restriction_{H_+}$, the \cstar\ $\Cal A_+$ and the
filtration  $\Cal F$.
{}From \cite{3, Theorem 4.5} we conclude that for every $f\in C_0(\Bbb R)$,
$$
\lim_{n\to\infty}{1\over
{d_n}}(f(\lambda_1^n)+f(\lambda_2^n)+\dots+f(\lambda_{d_n}^n))
= \tau_+(f(A_+)).
\leqno{2.9}
$$
{}From Propostion 2.1 and the fact that $f$ is continuous we have
$$
f(A_+)=P_+f(A)\!\restriction_{H_+} +\ \roman{compact}.
$$
Using the formula 2.6, we can rewrite the right side of 2.9 as follows
$$
\tau_+(P_+f(A)\!\restriction_{H_+} +\ \roman{compact})=\tau(f(A))=
\int_a^b f(x)\, d\mu_A(x),
$$
which establishes Theorem 2.3\qed
\enddemo

\subheading{3. Applications}

In order apply the general results of section 2
to establish the formula 1.4, we recall
some of the lore of almost periodic sequences
of complex numbers.  Consider the commutative \cstar\
$l^\infty =l^\infty(\Bbb Z)$.  The additive group $\Bbb Z$ acts naturally
on $l^\infty$ by translations, giving a $C^*$-dynamical
system.  An element $a\in l^\infty$ is {\sl almost periodic}
if the set of all of $\Bbb Z$-translates of $a$ is a relatively
norm-compact subset of $l^\infty$.  The set $AP(\Bbb Z)$ of
all almost periodic sequences is a unital \cstar -subalgebra
of $l^\infty$ which is invariant under the action of $\Bbb Z$.

The characters of $\Bbb Z$ are the sequences
$$
e_\lambda(n) = \lambda^n,\qquad n\in\Bbb Z
$$
where $\lambda$ is a complex number of absolute value $1$.
Every $e_\lambda$ is almost periodic, and in fact $AP(\Bbb Z)$
is precisely the closed linear span of all characters.  More
generally, if $S$ is any translation-invariant linear subspace of
$AP(\Bbb Z)$, then the {\sl spectrum} of $S$ is defined as the
following subset of the unit circle
$$
\sigma(S) = \{\lambda\in\Bbb T : e_\lambda\in S\}.
$$
$\sigma(S)$ is nonempty if $S\neq 0$, and more generally every
norm-closed translation-invariant subspace $S$ of almost periodic
sequences obeys uniform spectral synthesis in that
$$
S = \overline{\roman{span}}\,\{e_\lambda : \lambda\in\sigma(S)\}.\leqno{3.1}
$$
If, in addition, $S$ is a unital \cstar\ subalgebra of $AP(\Bbb Z)$ then
$\sigma(S)$ is a subgroup of the circle group $\Bbb T$.
Typically,  $\sigma(S)$ fails to be closed in the usual topology of $\Bbb T$;
indeed, every subset of $\Bbb T$ arises as the spectrum of some $S$
(so long as one is working with an underlying set theory that is
appropriate for functional analysis).
It follows easily from this summary that {\sl an almost periodic sequence}
$a$ {\sl is periodic iff} $\sigma(a)$ {\sl generates a finite subgroup
of} $\Bbb T$.

Finally, the norm-closed convex hull of the set of translates of
any almost periodic sequence
$a$ contains exactly one constant sequence $M(a)\I$, and the scalar
$M(a)$ is the von Neumann mean of $a$:
$$
M(a) = \lim_{n\to\infty}\frac{1}{2n+1}(a_{-n}+a_{-n+1}+\dots +a_n).
$$
Moreover, the convergence to $M(a)$ is {\sl uniform} in that
$$
\sup_{k\in\Bbb Z}|M(a) -\frac{1}{2n+1}(a_{k-n}+a_{k-n+1}+\dots +a_{k+n})| \to 0
$$
as $n\to\infty$.  $M$ defines a translation-invariant state of the
\cstar\ $AP(\Bbb Z)$.  We may conclude from these facts that if $C$ is any
unital translation-invariant $C^*$--subalgebra of $AP(\Bbb Z)$, then
the restriction of $M$ to $C$ is the {\sl only} translation-invariant
state on $C$.  A convenient reference for most
of the above is \cite{11}.

The following result is a consequence of standard techniques in
operator algebras, together with the preceding observations.
Since it is the essential link between the applications and
the material of section 2, we include a proof.

\proclaim{Proposition 3.2} Let $d=\{d_n: n\in\Bbb Z\}$ be a bounded
real sequence which is almost periodic but not periodic and
let $D$ be the corresponding diagonal operator acting on
$\l^2(\Bbb Z)$.  Then the \cstar\ generated by $D$ and the
bilateral shift is a simple \cstar\ possessing a unique tracial
state.
\endproclaim
\demo{proof}  Let $U$ be the bilateral shift acting on $\l^2(\Bbb Z)$,
$Uf(n) = f(n-1)$, $n\in\Bbb Z$.  For $f\in\l^2(\Bbb Z)$ we will
write $M_f$ for the obvious diagonal operator.  Let $\Cal A$ be the
\cstar\ generated by $M_d=D$ and $U$.  We first give a description
of $\Cal A$ which is more convenient for our purposes.

Let $\Delta$ be the translation-invariant $C^*$-subalgebra of
$AP(\Bbb Z)$ generated by $d$ and the constant sequence $1$, and
let $\Cal D$ be the set of all diagonal operators $M_f, f\in\Delta$.
Since $U^nM_fU^{-n}=M_{f_n}$ where $f_n$ is the translate of $f$ by $n$,
it follows that $\Cal A$ is the norm-closure of the set of
all operators which are finite sums of the form
$$
T = D_{-n}U^{-n}+D_{-n+1}U^{-n+1}+\dots+D_nU^n
$$
where each $D_k$ belongs to $\Cal D$ and $n=1,2,\dots$.
$\Cal A$ is therefore the image of the $C^*$-algebraic crossed product
$\Bbb Z \times \Delta$ under the obvious representation it has on
$l^2(\Bbb Z)$.  Thus, to show that $\Cal A$ is simple it
suffices to show that
$\Bbb Z\times \Delta$ is simple.  By the
characterization of simplicity of discrete crossed products
given in \cite{13, Theorem 8.11.12}, it is enough
to show that
\roster
\item
the only closed translation-invariant nonzero ideal in
$\Delta$ is $\Delta$ itself.
\item
The Connes spectrum of the action of $\Bbb Z$ on $\Delta$ is
all of $\Bbb T$.
\endroster
(1) follows immediately from the preceding remarks.  Indeed,
the preceding remarks imply that
any translation-invariant ideal must contain enough
unitary elements to span itself.  Hence if it is nonzero
it must be all of $\Delta$.  To prove (2),
note that since (1) is valid and
$\Delta$ is abelian, the Connes spectrum is identical
with the Arveson spectrum, and the latter is simply
the closure in $\Bbb T$ of the subgroup
$$
\sigma(\Delta)=\{\lambda\in\Bbb T: e_\lambda\in\Delta\}.
$$
But since $\Delta$ contains the sequence $d$ which is
{\sl not} periodic, $\sigma(\Delta)$ cannot be a finite
subgroup of $\Bbb T$; hence its closure is $\Bbb T$.

It remains to show that the crossed product $\Bbb Z\times \Delta$
has a unique tracial state.  Indeed, the translation-invariant
state of $\Delta$ gives rise to a tracial state $\tau$ of
$\Bbb Z\times \Delta$, and therefore of $\Cal A$.
On the other hand, suppose that $\rho$ is an arbitrary
tracial state on $\Cal A$.  For each $\lambda\in\sigma(\Delta)$,
consider the unitary operator $V_\lambda=M_{e_\lambda}$.  $V$ is
a unitary representation of the discrete abelian group $\sigma(\Delta)$,
$\Cal D$ is spanned by $\{V_\lambda:\lambda\in\sigma(\Delta)\}$,
and moreover $V$ satisfies the following commutation relation with $U$
$$
V_\lambda U = \lambda UV_\lambda , \qquad \lambda\in\sigma(\Delta).
$$
This implies that for every $n\in\Bbb Z$ and $D\in\Cal D$
we have
$$
\rho(DU^n) = \rho(V_\lambda DU^nV_{-\lambda}) = \lambda^n\rho(DU^n),
$$
for all $\lambda\in\sigma(\Delta)$.  If $n\neq0$ then it follows
that $\rho(DU^n)=0=\tau(DU^n)$.
If $n=0$, then the restriction of $\rho$ to
$\Cal D$ gives rise to a translation-invariant state of $\Delta$
and hence it agrees with the restriction of $\tau$ to $\Cal D$.
Hence $\rho=\tau$
\qed
\enddemo

We can now deduce applications to the problem of computing the
spectra of discretized Hamiltonians.  Let $d=\{d_n: n\in\Bbb Z\}$
be an almost periodic sequence of reals which is not periodic.
Consider the Hilbert space
$H=l^2(\Bbb Z)$, let $\{e_n: n\in\Bbb Z\}$ be the obvious basis
for $H$, and let $T$ be the tridiagonal operator defined by 1.2.
Let $\Cal A\subseteq\Cal B(H)$ be the \cstar\ generated by the
diagonal operator $M_d$ and the bilateral shift.  According to
3.2, $\Cal A$ is a simple \cstar\ with a unique tracial state
$\tau$.  $T$ certainly belongs to $\Cal A$, and we may consider
its spectral distribution $\mu_T$, defined by
$$
\int_{-\infty}^{+\infty} f(x)\,d\mu_T(x) = \tau(f(T)), \qquad f\in C_0(\Bbb R).
\leqno{3.3}
$$

\proclaim{Theorem 3.4} The spectrum of $T$ is the closed
support of $\mu_T$.  Moreover, for every positive integer $n$
let $\lambda_1^n < \lambda_2^n <\dots <\lambda_n^n$
be the eigenvalue list of the symmetric $n\times n$ matrix
$$
\pmatrix
d_{1}&1&0&\hdots&0&0\\
1&d_{2}&1&\hdots&0&0\\
0&1&d_{3}&\hdots&0&0\\
\vdots&\vdots&\vdots&\ddots&\vdots&\vdots\\
0&0&0&\hdots&d_{n-1}&1\\
0&0&0&\hdots&1&d_n
\endpmatrix.
\leqno{3.5}
$$
Then for every $f\in C_0(\Bbb R)$ we have
$$
\lim_{n\to\infty}\frac{1}{n}(f(\lambda_1^n)+f(\lambda_2^n)+\dots,
+f(\lambda_n^n)) = \int_{-\infty}^{+\infty} f(x)\,d\mu_T(x).
$$
\endproclaim
\demo{proof} Since $\Cal A$ is a simple \cstar, $\tau$ must
be a faithful trace.  Thus it is apparent from definition 3.3
that $\mu_T$ must be supported precisely on the spectrum of
$T$.

Let $\Cal F=\{H_1\subseteq H_2\subseteq\dots\}$ be the
filtration of $H$ defined by $H_n=[e_1,e_2,\dots,e_n]$.
Then $\Cal F$ is an improper filtration with respect to
which the bilateral shift has degree 1 and every diagonal
operator has degree 0.  Hence $\Cal F$ is an
$\Cal A$--filtration.  The hypothesis 2.4 of Theorem 2.3
is automatic in this case because of the simplicity of
$\Cal A$.  Thus Theorem 3.4 follows after an application of
Theorem 2.3 and its accompanying remarks\qed
\enddemo

\remark{Remarks} We reiterate some remarks from \cite{3} which
reveal the significance of Theorem 3.4 for computations.
For every $n=1,2,\dots$ and any Borel set $S\subseteq\Bbb R$,
let $N_n(S)$ be the number of eigenvalues of the matrix
2.5 which belong to $S$.  $N_n$ is a positive integer-valued
measure having total mass $n$.  Theorem 3.5 asserts that the
sequence of probability measures $n^{-1}N_n$ converges to
$\mu_T$ in the weak$^*$-topology of $C_0(\Bbb R)$.  Notice that
this is the same notion of convergence that one has in the
central limit theorem, and one may apply similar interpretations
in the context of this paper.

For example, suppose that $\lambda$ is a point in
the spectrum of $T$.  Then for every open neighborhood $I$
of $\lambda$ whose endpoints have $\mu_T$--measure zero,
and for any positive numbers $\alpha$, $\beta$ satisfying
$\alpha<\mu_T(I)<\beta$ we have
$$
n\alpha \leq N_n(I)\leq n\beta
$$
for sufficiently large $n$.  Notice that since $\alpha$ and
$\beta$ may be chosen arbitrarily close to $\mu_T(I)$,
we have a precise estimate of the rate of growth of
$N_n(I)$ with $n$.
If, on the other hand, $\lambda$ is not in the spectrum of
$T$, then for a sufficiently small interval $I$ containing
$\lambda$ we will have
$$
\lim_{n\to\infty}\frac{1}{n}N_n(I) = \mu_T(I)=0.
$$
Actually, in this case one can say more,
namely that the numbers $N_1(I), N_2(I),\dots$ are {\sl uniformly
bounded}.  Indeed, since $\Cal A$ is simple it cannot contain
any nonzero compact operators, and hence the spectrum of $T$ is
the same as its essential spectrum.  The assertion follows
from \cite{3, Theorem 3.8}.  In light of the
fact that there are extremely
fast algorithms for calculating the value of $N_n(S)$ for
large but fixed values of $n$ whenever
$S$ has the form $S=(-\infty,x]$ for $x\in\Bbb R$, these
remarks provide the basis for a very efficient method of
calculating the spectrum of $T$ \cite{4}.
\endremark

%
%
%
%
%\newpage

\enddocument